\documentclass[12pt]{amsart}
\usepackage{amssymb, amscd}

\newtheorem{Lemma}             {Lemma}
\newtheorem{Corollary}  [Lemma]{Corollary}
\newtheorem{Proposition}[Lemma]{Proposition}

\newtheorem{Example}    [Lemma]{Example}
\newtheorem{Theorem}    [Lemma]{Theorem}
\newtheorem{Conjecture} [Lemma]{Conjecture}
\newtheorem{Definition} [Lemma]{Definition}

\newcommand{\ov}{\overline}
\newcommand{\op}{\operatorname}

\title[Quadratic PIMs]{Quadratic principal indecomposable modules and strongly real elements of finite Groups}

\author{Rod Gow}
\address{School of Mathematics and Statistics\\University College Dublin\\IRELAND}
\email{Rod.Gow@ucd.ie}

\author{John C. Murray}
\address{Department of Mathematics and Statistics\\National University of Ireland Maynooth\\IRELAND}
\email{John.Murray@nuim.ie}

\date{September 29, 2017}

\begin{document}
\maketitle
\thispagestyle{empty}

\begin{abstract}
Let $P$ be a principal indecomposable module of a finite group $G$ in characteristic $2$ and let $\varphi$ be the Brauer character of the corresponding simple $G$-module. We show that $P$ affords a non-degenerate $G$-invariant quadratic form if and only if there are involutions $s,t\in G$ such that $st$ has odd order and $\varphi(st)/2$ is not an algebraic integer.

We then show that the number of isomorphism classes of quadratic principal indecomposable $G$-modules is equal to the number of strongly real conjugacy classes of odd order elements of $G$.
\end{abstract}

\section{Introduction}

Let $G$ be a finite group and let $p$ be a prime. Let ${\mathbb A}$ be the ring of algebraic integers in ${\mathbb C}$ and let ${\mathfrak M}$ be a maximal ideal of\/ ${\mathbb A}$ containing $p$. Set $R$ as the localisation of ${\mathbb A}$ at ${\mathfrak M}$. So $R$ is a complete discrete valuation ring with unique maximal ideal $J=R{\mathfrak M}$, $F:=\op{Frac}(R)$ is the field of algebraic numbers and $k:=R/J$ is the algebraic closure of $\op{GF}(p)$. Then $(F,R,k)$ is a so-called $p$-modular system for $G$. Note that $F$ and $k$ are splitting fields for all subgroups of $G$. We say that an $RG$-module has quadratic type if it affords a nondegenerate $G$-invariant quadratic form.

Let $P$ be a principal indecomposable $RG$-module. The principal indecomposable character of $P$ is the ordinary character of the $FG$-module $P\otimes_RF$. Now $P$ has a unique maximal submodule $\op{Soc}(P)$ containing $J(R)G$. The irreducible Brauer character corresponding to $P$ is the Brauer character of the $kG$-module $P/\op{Soc}(P)$. 

Recall that $g\in G$ is $p$-regular if its order is prime to $p$ and $p$-singular if its order is divisible by $p$. Also an element of $G$ is real if it is conjugate to its inverse, strongly real if it is inverted by an involution and weakly real if it is real but not strongly real. In particular the elements which square to the identity are strongly real.

Our main result is to give an efficient way of determining in characteristic $2$ if a self-dual principal indecomposable module has quadratic type:

\begin{Theorem}\label{T:strong+weak}
Let $p=2$ and let $P$ be a self-dual principal indecomposable $RG$-module with principal indecomposable character $\Phi$ and corresponding irreducible Brauer character $\varphi$. Then the following are equivalent:
\begin{itemize}
\item[(i)] $P$ has quadratic type.
\item[(ii)] $\varphi(g)\not\in2R$, for some strongly real $2$-regular $g\in G$.
\item [(iii)] $\frac{\Phi(g)}{|C_G(g)|}\in2R$, for all weakly real $2$-regular $g\in G$.
\end{itemize}
\end{Theorem}

As a consequence, we obtain:

\begin{Corollary}\label{T:strong_main}
If ${p=2}$ then the number of quadratic type principal indecomposable $RG$-modules equals the number of strongly real $2$-regular conjugacy classes of\/ $G$.
\end{Corollary}

Equivalently the number of non-quadratic type self-dual principal indecomposable $RG$-modules equals the number of weakly real $2$-regular conjugacy classes of $G$.

For $p$ odd, each principal indecomposable $RG$-module affords a non-degenerate $G$-invariant quadratic form or symplectic bilinear form. We do not know how to determine the number of each using the $p$-regular conjugacy classes of $G$.

\section{The type of principal indecomposable modules}

Prior work has classified the type of principal indecomposable $kG$-modules using character theory or ring theory, as we outline in this section. Let $(F,R,k)$ be as in the introduction. Then the group of units $\op{U}(R)$ of $R$ consists of all $p'$-roots of unity in ${\mathbb C}$. So the projection $R\rightarrow k$ induces an isomorphism $\op{U}(R)\cong k^\times$ of multiplicative groups (which of course depends on the choice ${\mathfrak M}$ of maximal ideal).

Let $M$ be an $RG$-module ($kG$-module) which is finitely generated and free as $R$-module. We say that $M$ has quadratic or symplectic type if $M$ affords a non-degenerate $G$-invariant $R$-valued ($k$-valued) quadratic form or symplectic bilinear form, respectively. It is known that when $p\ne2$ each indecomposable $kG$-module is either of quadratic type or of symplectic type. When $p=2$, each quadratic type $kG$-module is of symplectic type, but not conversely. 

The ring multiplication in $RG$ makes it into a module over itself, called the regular $RG$-module. The direct summands of this module are called principal indecomposable $RG$-modules. We say that a principal indecomposable $RG$-module $P$ is trivial if $P/\op{Rad}(P)$ is the trivial $kG$-module. As was shown by R. Brauer, the number of principal indecomposable $RG$-modules equals the number of $p$-regular (elements of order prime to $p$) conjugacy classes of $G$ and the number of self-dual principal indecomposable $RG$-modules equals the number of real $p$-regular conjugacy classes of $G$.

Each self-dual irreducible $FG$-module is orthogonal or symplectic as it affords a non-degenerate $G$-invariant $F$-valued symmetric bilinear form or symplectic bilinear form, respectively. We say that the corresponding irreducible character of $G$ has orthogonal or symplectic type, respectively. The type can be detected by computing the Frobenius-Schur (F-S) indicator of the character; symmetric type irreducible characters have F-S indicator $+1$ and symplectic type irreducible characters have F-S indicator $-1$.

Suppose that $P$ is a self-dual principal indecomposable $RG$-module and let $\Phi$ be the principal indecomposable character of $P$. Then $\Phi$ is real-valued. Suppose first that $p$ is odd. Then by \cite{Wi76} and \cite{Th84} some real-valued irreducible character of $G$ occurs with odd multiplicity in $\Phi$. Moreover $P$ has quadratic or symplectic type, as this character has orthogonal or symplectic type, respectively.

Suppose then that $p=2$. Then by \cite{GW93}, some orthogonal irreducible character of $G$ occurs with odd multiplicity in $\Phi$.  As a consequence, $P$ may have quadratic type but it cannot have symplectic type. If some symplectic irreducible character of $G$ occurs with odd multiplicity in $\Phi$, then $P$ cannot have quadratic type. However $P$ may be of non-quadratic type without the occurrence of such a character.

Reduction mod $J$ induces a bijection between the principal indecomposable $RG$-modules and the principal indecomposable $kG$-modules. Under this bijection a principal indecomposable $RG$-module has quadratic type if and only if the corresponding principal indecomposable $kG$-module has quadratic type.  So for convenience we can and do work with principal indecomposable $kG$-modules.

Write $1_G=e_1+\dots+e_m$ where the $e_i$ are pairwise orthogonal primitive idempotents in $kG$ (i.e. $e_ie_j=\delta_{ij}e_i$, for $1\leq i,j\leq m$ and $m\geq1$ is as large as possible). Then each $kGe_i$ is a principal indecomposable $kG$-module, and all principal indecomposable $kG$-modules have this form. Now the map $g\rightarrow g^{-1}$, for $g\in G$, extends to an involutary $k$-algebra anti-automorphism $^o$ of $kG$, called the contragredient map. Each $e_i^o$ is a primitive idempotent in $kG$ and $kGe_i^o$ is isomorphic to the dual module $(kGe_i)^*:=\op{Hom}(kGe_i,k)$ of $kGe_i$.

We must distinguish between $p$ odd and $p=2$:

\begin{Lemma}[Landrock-Manz]\label{L:Landrock-Manz}
Let $p$ be an odd prime and let $e$ be a primitive idempotent in $kG$. Then $kGe$ has quadratic type if and only if there is a primitive idempotent $f\in kG$, with $kGf=kGe$, such that $f=f^o$.
\end{Lemma}

In contrast, if $p=2$, Gow and Willems showed that there is a primitive idempotent $f\in kG$,  with $kGf=kGe$, such that $f=f^o$ if and only if $kGe$ is the projective cover of the trivial $kG$-module. Moreover, they proved the following analogue of Lemma \ref{L:Landrock-Manz}:

\begin{Lemma}[Gow-Willems]\label{L:tf^ot=f}
Let $p=2$ and let $e$ be a primitive idempotent in $kG$. Then $kGe$ has quadratic type if and only if there is an involution $t\in G$ and a primitive idempotent $f\in kG$, with $kGf=kGe$, such that $t^{-1}ft=f^o$.
\begin{proof}[Outline Proof]
We may assume that $kGe$ is not the projective cover of the trivial $kG$-module. Then it is known that each $G$-invariant symmetric bilinear form on $kGe$ is symplectic and is the polarization of a $G$-invariant quadratic form on $kGe$ (see \cite[Proposition 2.2]{GW93}).

Let $B_1$ be the symmetric bilinear form on $kG$ with respect to which the elements of $G$ form an orthonormal basis of $kG$. Then $B_1$ is non-degenerate and $G$-invariant. Its adjoint is the contragredient map. Next recall that $\op{End}_{kG}(kG)$ can be identified with the opposite ring $kG^{op}$. Here $x\in kG^{op}$ defines the $kG$-homomorphism $y\rightarrow yx$, for all $y\in kG$. So if we define $B_x(y,z):=B_1(yx,z)$, for all $y,z\in kG$, then $B_x$ is a $G$-invariant bilinear form on $kG$ and $\{B_x\mid x\in kG\}$ give all $G$-invariant bilinear forms on $kG$. Moreover, $B_x$ is non-degenerate, symmetric or symplectic as $x$ is a unit in $kG$, $x=x^o$ or $x=x^o$ and $x_1=0$, respectively.

Next $(kGe)^*\cong kGe^o$. So the space of $G$-invariant bilinear forms on $kGe$ can be identified with $ekGe^o$. As a consequence, each $G$-invariant bilinear form on $kGe$ is the restriction  from $kG$ to $kGe$ of a $G$-invariant bilinear form $B_{exe^o}$, for some unique $exe^o\in ekGe^o$. Now suppose that $B_{exe^o}$ is a non-degenerate symplectic bilinear form on $kGe$. Write $exe^o=\sum_{g\in G}x_gg$. Then $x_1=0_k$ and $x_g=x_{g^{-1}}$, for all $g\in G$. As $kGe$ is indecomposable, $B_{x_g(g+g^{-1})}$ is degenerate on $kGe$. So there is an involution $t\in G$ such that $B_{x_tt}$ restricts to a non-degenerate (symplectic) bilinear form on $kGe$. Thus $x_t\ne0_k$ and $B_t$ restricts to a non-degenerate bilinear form on $kGe$. Let $f\in kG^{op}$ be the projection onto $kGe$, with kernel the complement of $kGe$ with respect to $B_t$. Now $B_t$ has adjoint $x\mapsto tx^ot$ on $kG^{op}$. So $f$ is a primitive idempotent in $kG$ with $kGe=kGf$ such that $tf^ot=f$.

Conversely, given an involution $t\in G$ and an idempotent $f\in kG$ as in the statement, it is easy to see that $B_t$ restricts to a non-degenerate symplectic bilinear form on $kGe$.
\end{proof}
\end{Lemma}

If $kGe$ has quadratic type and $t$ is as in the conclusion of the Lemma, we say that $t$ inverts $e$ or $kGe$, or the projective indecomposable character of $kGe$.

Fong's Lemma states that when $\op{char}(k)=2$ every non-trivial self-dual irreducible $kG$-module affords a non-degenerate $G$-invariant symplectic form. It is clear that this form is unique up to a non-zero scalar. In \cite[Proposition 5.8]{M16} the second author showed:

\begin{Lemma}\label{L:Murray}
Suppose that $p=2$ and that $M$ is a non-trivial self-dual irreducible $kG$-module, with Fong form $B$. Then the projective cover of $M$ has quadratic type if and only if $B(tm,m)\ne0$ for some $m\in M$ and involution $t\in G$.
\begin{proof}[Outline Proof]
The annihilator of $M$ in $kG$ is a maximal $2$-sided ideal $\op{Ann}(M)$ of $kG$ with $kG/\op{Ann}(M)\cong\op{End}_k(M)$. Suppose first that the principal indecomposable module of $M$ is orthogonal. Then by Lemma \ref{L:tf^ot=f} there is an involution $t\in G$ and a primitive idempotent $f\in kG$, such that $kGf/\op{Rad}(kGf)\cong M$ and $t^{-1}ft=f^o$. Set $x^{ot}:=tx^ot$, for all $x\in kG$. Then $^{ot}$ is an involutary $k$-algebra anti-automorphism of $kG$. Now $\op{Ann}(M)$ is invariant under $^{ot}$. So $^{ot}$ induces an involutary $k$-algebra anti-automorphism on $\op{End}_k(M)$. It is readily established that $g^{ot}=g^{-t}$, for all $g\in G$ in $\op{End}_k(M)$.

Next define $B_t(x,y):=B(tx,y)$, for all $x,y\in\op{End}_k(M)$. Then $B_t$ is a non-degenerate $\op{C}_G(t)$-invariant symmetric bilinear form on $M$. Its adjoint coincides with $^{ot}$ on $\op{End}_k(M)$, as by irreducibility of $M$, the image of $G$ spans $\op{End}_k(M)$. Now $f+\op{Ann}(M)$ is a primitive idempotent in $\op{End}_k(M)$ which is $^{ot}$-invariant. It follows that $B_t$ is non-degenerate on the $1$-dimensional subspace $fM$ of $M$ i.e. there is $m\in fM$ such that $B(tm,m)\ne0_k$.

Conversely, suppose that $B(tm,m)\ne0$ for some involution $t\in G$ and some $m\in M$. Let $\hat f$ be orthogonal projection onto $km$ with respect to the non-degenerate symmetric form $B_t$ on $M$. By idempotent lifting (c.f. \cite[Proposition 1.4]{LM}) there is a primitive idempotent $f\in kG$ such that $f^{ot}=f$ and $f+\op{Ann}(M)=\hat f$. Now $kGf/\op{Rad}(kGf)\cong M$. So $kGf$ is orthogonal, according to Lemma \ref{L:tf^ot=f}.
\end{proof}
\end{Lemma}

\section{Type of principal indecomposable modules of ${\mathbb R}$-elementary subgroups}

Recall that a group is elementary if it has the form $C\times P$, where $P$ is a $p$-group, for some prime $p$ and $C$ is a cyclic $p'$-group. Brauer's induction theorem states that every ${\mathbb C}$-character of $G$ is an integer combination of characters induced from linear characters of elementary subgroups of $G$.
For the rest of this section $p=2$.

Now an ${\mathbb R}$-elementary group has the form $C\rtimes P$, where $C$ is cyclic of odd order and $P$ is a $2$-group such that every element of $P$ centralizes or inverts $C$. The real version of Brauer's induction theorem (which is a special case of the Witt-Berman theorem) is that every ${\mathbb R}$-character of $G$ is an integer combination of characters induced from ${\mathbb R}$-characters of ${\mathbb R}$-elementary subgroups of $G$. In this section we determine the type of the principal indecomposable modules of ${\mathbb R}$-elementary groups.

As usual the centralizer $\op{C}_G(g)$ or extended centralizer $\op{C}_G^*(g)$ of $g\in G$ is the normalizer of $\{g\}$ or $\{g,g^{-1}\}$ in $G$, respectively. Then $\op{C}_G(g)\subseteq\op{C}_G^*(g)$ and $|\op{C}_G^*(g):\op{C}_G(g)|\leq2$. Moreover $g$ is strongly real if $g=g^{-1}$ or $\op{C}_G^*(g)\backslash\op{C}_G(g)$ contains an involution and weakly real if $\op{C}_G^*(g)$ does not split over $\op{C}_G(g)$. 

Our main result here is a special case of a result on the principal indecomposable modules of solvable groups which is due to G. Navarro and the second author \cite[Theorem 15]{MN16}:

\begin{Proposition}\label{P:E/Dsplits}
Let $g\in G$ and let $E\in\op{Syl}_2(\op{C}_G^*(g))$. Then a non-trivial principal indecomposable $R(\langle g\rangle\rtimes E)$-module has quadratic type if and only if\/ $g$ is strongly real.
\begin{proof}
We may assume that $g$ is real in $G$. Set $C:=\langle g\rangle$ and $H:=C\rtimes E$. Let $P$ be a non-trivial principal indecomposable $RH$-module, and let $\Phi$ be the principal indecomposable character of $P$. Then $P\cong M{\uparrow^H}$, for some non-trivial $1$-dimensional $kC$-module $M$.

Suppose first that $g$ is strongly real. Let $t$ be an involution in $E$ which inverts $g$. Then $M{\uparrow^{\langle g,t\rangle}}$ is a self-dual irreducible $k(\langle g,t\rangle)$-module. Fong's Lemma implies that $M{\uparrow^{\langle g,t\rangle}}$ affords a symplectic geometry. Then the induced module $P\cong(M{\uparrow^{\langle g,t\rangle}}){\uparrow^H}$ affords the induced form. So $P$ has symplectic, hence quadratic type.

Conversely suppose that $P$ has quadratic type. Set $D=\op{C}_E(g)$ and $\hat M:=(\op{Inf}_{C\times D/D}^{C\times D}M){\uparrow^H}$. Then $\hat M$ is a self-dual irreducible $kH$-module which is isomorphic to $P/\op{Rad}(P)$. Let $B$ be a Fong form on $\hat M$. According to Lemma \ref{L:Murray}, $B(tm,m)\ne0_k$, for some involution $t\in H$ and some $m\in\hat M$. Now $\hat M{\downarrow_{C\times D}}=M_1\oplus M_2$, where $M_1$ is an irreducible $k(C\times D)$-module and $M_2\cong M_1^*\not\cong M_1$. Write $m=m_1+m_2$, where $m_1\in M_1$ and $m_2\in M_2$. Then
$$
\begin{aligned}
B(tm,m)
&=B(tm_1,m_1)+B(tm_1,m_2)+B(tm_2,m_1)+B(tm_2,m_2)\\
&=B(tm_1,m_1)+B(tm_2,m_2),\quad\mbox{as $B(tm_1,m_2)=B(tm_2,m_1)$.}
\end{aligned}
$$
So we may assume without loss of generality that $B(tm_1,m_1)\ne0_k$. As $M_1\not\cong M_1^*$, $B$ is identically zero on $M_1$. So $tm_1\not\in M_1$, which forces $t\in H\backslash(C\times D)$. Then $g^t=g^{-1}$. So $g$ is strongly real.
\end{proof}
\end{Proposition}

\section{Values of principal indecomposable and Brauer characters}

In this section $p=2$. We use Proposition \ref{P:E/Dsplits} to clarify the relationship between the strongly and weakly real $2$-regular conjugacy classes of\/ $G$ and the quadratic and non-quadratic self-dual principal indecomposable $RG$-modules.

Each $kG$-module $M$ has a Brauer character $\varphi_M$; if $g\in G$ has odd order then $M$ has a basis of eigenvectors of $g$. Each eigenvalue of $g$ on $M$ is a $2'$-roots of unity in $k$. So they can be lifted to $2'$-roots of unity in $R$, via the isomorphism $\op{U}(R)\cong k^\times$. Then $\varphi_M(g)$ is defined to be the sum of these roots of unity. In particular $\varphi_M$ is an ${\mathbb A}$-valued function defined on the $2$-regular elements of\/ $G$. As is standard, we use $\op{IBr}(G)$ to denote the Brauer characters of the irreducible $kG$-modules.

The restriction $\chi^*$ of an ordinary character $\chi$ of $G$ to the $p$-regular elements of\/ $G$ is a Brauer character. So $\chi^*=\sum_{\varphi\in\op{IBr}(G)}d_{\chi\varphi}\varphi$, where the $d_{\chi\varphi}$ are non-negative integers, called decomposition numbers. For $\varphi\in\op{IBr}(G)$ set $\Phi_\varphi:=\sum_{\chi\in\op{Irr}(G)}d_{\chi\varphi}\chi$. Then $\Phi_\varphi$ is the principal indecomposable character of\/ $G$ corresponding to $\varphi$. It is easy to see that $\Phi_\varphi$ vanishes on all $p$-singular elements of $G$.

Recall that a vertex of an indecomposable $kG$-module $M$ is a subgroup $V$ of\/ $G$ which is minimal subject to $M$ being a direct summand of a module induced from $V$ to $G$. J.~A.~Green showed that $V$ is a $p$-subgroup of\/ $G$, and moreover $V$ is uniquely determined up to $G$-conjugacy. The next result and its corollary are due to Green.

\begin{Lemma}
Let $M$ be an indecomposable $kG$-module, let $g\in G$ be $p$-regular and let $D\in\op{Syl}_p(\op{C}_G(g))$. Then $M$ has a vertex $V$ such that $\frac{\varphi_M(g)}{|D:V\cap D|}$ is an algebraic integer.
\end{Lemma}

In particular if\/ $\varphi_M(g)\not\in 2{\mathbb A}$ then $D$ is contained in some vertex of\/ $M$. 

\begin{Corollary}
Let $\Phi$ be a principal indecomposable character of\/ $G$ and let $g\in G$. Then $\frac{\Phi(g)}{|\op{C}_G(g)|_p}$ is an algebraic integer.
\begin{proof}
As $\Phi$ is zero on $p$-singular elements, we may assume that $g$ is $p$-regular. Now $\Phi^*$ is the Brauer character of a principal indecomposable $kG$-module and the trivial group is a vertex of this module. So $|\op{C}_G(g)|_p$ divides $\Phi(g)$ in ${\mathbb A}$.
\end{proof}
\end{Corollary}

Let $\Phi_1$ be the trivial principal indecomposable character of $G$ i.e. $\Phi_1$ corresponds to the trivial $2$-Brauer character $\varphi_1$ of $G$. In \cite{GW93} Gow and Willems proved that $\frac{\Phi_1(1)}{|G|_2}$ is odd. We complement this result with:

\begin{Theorem}\label{T:trivialPIM}
If\/ $g\in G$ is real and non-trivial then $\frac{\Phi_1(g)}{|\op{C}_G(g)|_2}$ is twice an algebraic integer.
\begin{proof}
By the second orthogonality relation $0=\hspace{-.25cm}\sum\limits_{\chi\in\op{Irr}(G)}\chi(1)\chi(g)=\hspace{-.25cm}\sum\limits_{\varphi\in\op{IBr}(G)}\varphi(1)\Phi_\varphi(g)$. So
$$
\frac{\Phi_1(g)}{|\op{C}_G(g)|_2}=-\sum_{\varphi\in\op{IBr}(G)\atop\varphi\ne\varphi_1}\varphi(1)\frac{\Phi_\varphi(g)}{|\op{C}_G(g)|_2}.
$$
For $\varphi\ne\ov\varphi$ we get two equal summands. So their sum cancels mod $2$. If\/ $\varphi=\ov\varphi$ and $\varphi\ne\varphi_1$, Fong's Lemma implies that $\varphi(1)$ is even. The conclusion follows from these facts.
\end{proof}
\end{Theorem}

\begin{Lemma}\label{L:weakly_real}
Let $g$ be a weakly real element of\/ $G$ and let $\Phi$ be a quadratic type principal indecomposable character of\/ $G$. Then $\frac{\Phi(g)}{|\op{C}_G(g)|_2}\in2{\mathbb A}$. 
\begin{proof}
We may assume that $g$ is $2$-regular. Suppose for the sake of contradiction that $\frac{\Phi(g)}{|\op{C}_G(g)|_2}\not\in2{\mathbb A}$. Let $E\in\op{Syl}_2(\op{C}_G^*(g))$ and set $H:=\langle g\rangle\rtimes E$. Write $\Phi{\downarrow_H}=\sum\alpha_\Lambda\Lambda$, where $\Lambda$ ranges over the principal indecomposable characters of\/ $H$. Lemma \ref{T:trivialPIM} implies that there is a non-trivial principal indecomposable character $\Lambda$ such that $\alpha_\Lambda$ is odd and $\frac{\Lambda(g)}{|\op{C}_E(g)|}\not\in2{\mathbb A}$.
 
Let $P$ be the projective indecomposable $kG$-module corresponding to $\Phi$ and let $Q$ be the principal indecomposable $kH$-module corresponding to $\Lambda$. Then $Q$ is a self-dual module which occurs with odd multiplicity $\alpha_\Lambda$ as a direct summand of $P{\downarrow_H}$. As $P$ is orthogonal, so too is $Q$. So $g$ is strongly real, according to Proposition \ref{P:E/Dsplits}. This contradiction completes the proof.
\end{proof}
\end{Lemma}

We can refine this result, using the techniques developed in \cite{M16}:

\begin{Lemma}\label{L:inverts}
Let $g$ be a real element of $G$ and let $\Phi$ be a quadratic principal indecomposable character of\/ $G$ such that $\frac{\Phi(g)}{|\op{C}_G(g)|_2}\not\in2{\mathbb A}$. Suppose that $t\in G$ is an involution which inverts $\Phi$ (c.f. the comment after Lemma \ref{L:tf^ot=f}). Then some conjugate of\/ $t$ inverts $g$.
\begin{proof}[Sketch Proof]
By Lemma \ref{L:tf^ot=f} there is a primitive idempotent $f$ in $kG$ such that $kGf$ is the principal indecomposable $kG$-module corresponding to $\Phi$ and $B_t$ restricts to a non-degenerate symplectic form on $kGf$. Let $H=\langle g\rangle\rtimes E$ and $Q$ be as in Lemma \ref{L:weakly_real}; so $Q$ is a projective indecomposable $kH$-module which occurs with odd multiplicity in $kGf{\downarrow_H}$. Then $Q$ is a non-degenerate component of the symplectic module $(kGf,B_t){\downarrow_H}$. Now $(kG,B_t)=(k\langle t\rangle,B_t){\uparrow^G}$. So $Q$ is a non-degenerate component of
\begin{equation}\label{E:mackey}
(k\langle t\rangle,B_t){\uparrow^G}{\downarrow_H}=\mathop{\perp}\limits_{\langle t\rangle aH\subseteq G}(k\langle t^a\rangle,B_{t^a}){\downarrow_{\langle t^a\rangle\cap H}}{\uparrow^H}.
\end{equation}
Let ${a\in G}$. First suppose that $\langle t^a\rangle\subseteq H$. Then the right hand side has a summand
\begin{equation}\label{E:tinH}
{(k\langle t^a\rangle,B_{t^a}){\uparrow^H}\!\cong\!(kH,B_{t^a})}.
\end{equation}
The other possibility is that $t^a\not\in H$. Then $\langle t^a\rangle\cap H=1$ and $(k\langle t^a\rangle,B_{t^a}){\downarrow_1}\cong(k^2,\hat{B_1})$ is a symplectic plane. Let $1=e_1+\dots+e_u$ be a decomposition of\/ $1$ as a sum of pairwise orthogonal primitive idempotents in $kH$. For each $i$, the bilinear form $B_1$ on $kG$ defines a perfect $G$-pairing $kGe_i\times kGe_i^o\rightarrow k$. We set $(kGe_i\oplus kGe_i^o,\hat{B_1})$ as the corresponding symplectic {\em paired} module. It is important to note that this module has no proper non-degenerate component. Then the summand on the right hand side of \eqref{E:mackey} has the orthogonal decomposition
\begin{equation}\label{E:tnotinH}
(k^2,\hat{B_1}){\uparrow^H}=\mathop{\perp}\limits_{i=1}^u(kHe_i\oplus kHe_i^o,\hat{B_1}).
\end{equation}
Now $Q$ is an indecomposable non-degenerate submodule of an orthogonal sum of modules of the form \eqref{E:tinH} and \eqref{E:tnotinH}. As the modules in \eqref{E:tnotinH} have no such submodules, we deduce that there is $a\in G$ such that $s:=t^a\in H$ and $Q$ is a non-degenerate submodule of\/ $(kH,B_s)$. Equivalently there is a primitive idempotent $e\in kH^{op}$ such that $kHe\cong Q$ and $e^{os}=e$. Now $k(C\times D)e$ is a projective indecomposable $k(C\times D)$-module which is not self-dual (as it is not the projective cover of the trivial $k(C\times D)$-module). So $s\not\in(C\times D)$. We conclude that $g^s=g^{-1}$.
\end{proof}
\end{Lemma}

Complementing Lemma \ref{L:weakly_real}, we have:

\begin{Lemma}\label{L:strongly_real}
Let $g$ be a strongly real $2$-regular element of\/ $G$ and let $P$ be a non-quadratic self-dual principal indecomposable $RG$-module with associated ${\varphi\in\op{IBr}(G)}$. Then ${\varphi(g)\in2{\mathbb A}}$.
\begin{proof}
We note that $\varphi(1)$ is even by Fong's Lemma. So we may assume that $g\ne1$. As before let $E\in\op{Syl}_2(\op{C}_G^*(g))$ and set $H:=\langle g\rangle\rtimes E$.
 
Suppose for the sake of contradiction that $\varphi(g)\not\in 2{\mathbb A}$. Write $\varphi{\downarrow_H}=\sum\beta_\lambda\lambda$, where $\lambda$ ranges over the irreducible Brauer characters of\/ $H$. So by hypothesis there is $\lambda\in\op{IBr}(H)$ such that $\beta_\lambda$ is odd and $\lambda(g)\not\in 2{\mathbb A}$.
 
Let $Q$ be the projective indecomposable $kH$-module corresponding to $\lambda$. As $g$ is strongly real, it follows from Proposition \ref{P:E/Dsplits} that $Q$ has quadratic type. Now by Frobenius-Nakayama reciprocity $P$ occurs with odd multiplicity $\beta_\lambda$ in $Q{\uparrow^G}$. As $P\cong P^*$, we conclude that $P$ has quadratic type. This contradiction completes the proof.
\end{proof}
\end{Lemma}

Let $M$ be the irreducible $kG$-module whose Brauer character is $\varphi$. In \cite{M16} the second author assigned a symplectic vertex $T$ to $M$; $T$ is a minimal subgroup of $G$ such that $M$ is an orthogonal direct summand of a symplectic $kT$-module induced up to $G$. Just as with Green vertices, $T$ is uniquely determined up to $G$-conjugacy. In view of Lemma \ref{L:inverts} we hazard the following, which would complement Lemma \ref{L:inverts} and strengthen Lemma \ref{L:strongly_real}:

\begin{Conjecture}
Let $M$ be the irreducible $kG$-module, let $\varphi$ be the Brauer character of $M$ and let $g$ be a real $2$-regular element of\/ $G$ such that $\varphi(g)\not\in2{\mathbb A}$. Then $M$ has an extended defect group which contains an extended defect group of $g$ but which is not contained in $\op{C}_G(g)$.
\end{Conjecture}

\section{Type of principal indecomposable $RG$-modules}

Recall that $R$ is the localisation of the ring of algebraic integers at a maximal ideal containing $2$. In particular $R$ has a unique maximal ideal $J$. Suppose that $G$ has $\ell$\, $2$-regular conjugacy classes, $r$ real $2$-regular conjugacy classes and $s$ strongly real $2$-regular conjugacy classes. So $G$ has $\ell$ irreducible $2$-Brauer characters and $r$ real irreducible $2$-Brauer characters. Suppose also that $G$ has $\sigma$ quadratic type principal indecomposable $RG$-modules.

List the irreducible $2$-Brauer characters of $G$ as $\varphi_1,\dots,\varphi_\ell$ and let $g_1,\dots,g_\ell$ be a set of representatives for the $2$-regular conjugacy classes of $G$. Set $\Phi_k$ as the principal indecomposable character of $G$ corresponding to $\varphi_k$, for $k=1,\dots,\ell$. 

The second orthogonality relation (see Theorem \ref{T:trivialPIM}) give equations in $R$:
$$
\sum_{k=1}^\ell\frac{\Phi_k( g_i^{-1} )}{|\op{C}_G(g_i)|}\varphi_k(g_j)=\delta_{ij},\quad\mbox{for $1\leq i,j\leq\ell$}.
$$
Suppose that $g_i$ and $g_j$ are real in $G$. Then $\frac{\Phi_k( g_i^{-1} )}{|\op{C}_G(g_i)|}\varphi_k(g_j)=\frac{\ov\Phi_k( g_i^{-1} )}{|\op{C}_G(g_i)|}\ov\varphi_k(g_j)$. So the contribution of the non-real $\varphi_k$ to the above displayed sum is zero, modulo $J$.

We may choose our notation so that $\varphi_1,\dots,\varphi_r$ are the real irreducible $2$-Brauer characters of $G$ and $g_1,\dots,g_r$ are in the real $2$-regular conjugacy classes of $G$. Define the $r\times r$-matrices:
$$
A=\left[\frac{\Phi_j(g_i^{-1})}{|\op{C}_G(g_i)|}\right],\quad
B=\left[\varphi_i(g_j)\right].
$$
They involve only the real principal indecomposable characters, the real irreducible Brauer characters and the real $2$-regular elements of\/ $G$. By the work above $AB\equiv I$ (mod $J$).

We further refine our notation so that $\varphi_1,\dots,\varphi_\sigma$ are the Brauer characters of the quadratic principal indecomposable $RG$-modules and $g_1,\dots,g_s$ are in the strongly real $2$-regular conjugacy classes of $G$. Thus $\varphi_{\sigma+1},\dots,\varphi_r$  are the Brauer characters of the non-quadratic principal indecomposable $RG$-modules and $g_{s+1},\dots,g_r$ are in the weakly real $2$-regular conjugacy classes of $G$. The matrices $A$ and $B$ have corresponding block forms:
$$
A=\left[\begin{array}{ll}A_{11}&A_{12}\\A_{21}&A_{22}\end{array}\right],\quad
B=\left[\begin{array}{ll}B_{11}&B_{12}\\B_{21}&B_{22}\end{array}\right].
$$
So $A_{11}$ is the $s\times\sigma$ submatrix of\/ $A$ with rows and columns indexed by the strongly real classes and quadratic type principal indecomposable characters. Likewise $B_{11}$ is the $\sigma\times s$ submatrix of $B$ with rows and columns indexed by Brauer characters of quadratic principal indecomposable $RG$-modules and strongly real classes, respectively.

Now each term in $A_{21}$ has the form $\frac{\Phi_k(g_i^{-1})}{|\op{C}_G(g_i)|}$ where $g_i$ is weakly real and $\Phi_k$ has quadratic type. Likewise each term in $B_{21}$ has the form $\varphi_k(g_j)$ where $\varphi_k$ is the Brauer character of a non-quadratic principal indecomposable $RG$-module and $g_j$ is strongly real. So according to Lemmas \ref{L:weakly_real} and \ref{L:strongly_real}, all these terms belong to $J$. Thus
$$
A_{21}\equiv0\,\mbox{(mod $J$)}\quad\mbox{and}\quad
B_{21}\equiv0\,\mbox{(mod $J$)}.
$$
It follows from this that $A_{11}B_{11}\equiv I$ (mod $J$) and $A_{22}B_{22}\equiv I$ (mod $J$). In particular $A_{11}$ and $A_{22}$ have full row rank (mod $J$). So $s\leq\sigma$ and $r-s\leq r-\sigma$. We conclude that $s=\sigma$. This proves our main theorem, which we restate here for the convenience of the reader:

\begin{Theorem}\label{T:=}
The number of strongly real $2$-regular conjugacy classes of\/ $G$ equals the number of quadratic type principal indecomposable $RG$-modules and the number of weakly real $2$-regular conjugacy classes of\/ $G$ equals the number of non-quadratic type self-dual principal indecomposable $RG$-modules.
\end{Theorem}

The analysis above translates into the following criteria for strong and weak reality, which proves Theorem \ref{T:strong+weak}:

\begin{Proposition}\label{P:strong+weak}
Let $p=2$ and let $P$ be a self-dual principal indecomposable $RG$-module with principal indecomposable character $\Phi$ and corresponding irreducible Brauer character $\varphi$. Then:
\begin{itemize}
\item[(i)] $P$ has quadratic type if and only if $\varphi(g)\not\in2R$, for some strongly real $2$-regular $g\in G$.
\item[(ii)] $P$ has non-quadratic type if and only if $\frac{\Phi(g)}{|C_G(g)|}\not\in2R$, for some weakly real $g\in G$.
\end{itemize}
\end{Proposition}

\section{Strong and weak projective indecomposable modules}

In this section we give examples of strong and weak principal indecomposable characters, many of which involve quasisimple finite groups. In addition to \cite{Atlas}, we use the notation and decomposition matrices provided by the Modular Atlas homepage \cite{MOC}.

Let $g\in G$ be $2$-regular. If $g$ is strongly real, $tgt=g^{-1}$, for some involution $t\in G$. Then $s:=gt$ is an involution which is conjugate to $t$ in the dihedral group $\langle g,t\rangle$ and $g=st$. Conversely, if $g=uv$ where $u,v\in G$ are involutions, then $v$ inverts $g$ and so $g$ is strongly real. So by the class multiplication formula, $g$ is strongly real if and only if
$$
\sum_{\chi\in\op{Irr}(G)}\frac{\chi(t)^2\chi(g)}{\chi(1)}\ne0,
$$
for some involution $t\in G$. We note that the character table of $G$ determines the prime divisors of the orders of the elements of $G$. In particular the character table determines the $2$-regular conjugacy classes of $G$. However the character table of a group does not generally determine which conjugacy classes of $2$-elements are involutions e.g. $D_8$ and $Q_8$ have the `same' character tables.

\begin{Example}
$G=2.A_5$. All three non-trivial $2$-regular classes of $G$ are weakly real, for example because $G$ has a unique involution. So all three non-trivial $2$-principal indecomposable characters $\Phi_2,\Phi_3$ and $\Phi_4$ of $2.A_5$ are of non-quadratic type. This also follows from the fact that $\Phi_2$ contains the  symplectic irreducible characters $\chi_6$ and $\chi_9$ with odd multiplicity $1$, $\Phi_3$ contains the  symplectic irreducible characters $\chi_7$ and $\chi_9$ with odd multiplicity $1$ and $\Phi_4$ contains the  symplectic irreducible character $\chi_8$ with odd multiplicity $1$.
\end{Example}

\begin{Example}
$G=\op{Sp}(4,5)$, or $2.\op{S}_4(5)$ in \cite{Atlas} notation. Then $G$ has $60$ irreducible ${\mathbb C}$-characters, $26$ of which are faithful. All non-faithful characters are orthogonal and all faithful characters are symplectic. Now $G$ has fifteen $2$-regular classes, all of whom are real. By examination of the $2$-decomposition matrix, each of the $2$-principal indecomposable characters $\Phi_i$, for $i\in\{2,3,4,6,8,9,10,11,12,13,14\}$ contain at least one symplectic irreducible constituent with odd multiplicity. So these twelve principal indecomposable characters are non-quadratic. By computing the square of the involution class $2B$, we see that $1A$, $3A$ and $5D$ are the strongly real $2$-regular classes of $G$. Now $\varphi_5=\chi_6^*-\chi_1^*$ and $\chi_6(5D)=0$. So $\varphi_5(5D)=-1$ is odd. Similarly $\varphi_7=\chi_{11}^*$ and $\chi_{11}(3A)=5$ is odd. So $\Phi_1,\Phi_5$ and $\Phi_7$ are the quadratic principal indecomposable characters of $G$.
\end{Example}

\begin{Example}
$G=McL$ has $24$ irreducible ${\mathbb C}$-characters, of whom two, $\chi_{11}$ and $\chi_{13}$, are symplectic and $10$ are orthogonal. Also $G$ has $11$ $2$-regular conjugacy classes, $5$ of whom are real. Examination of the square of the involution class $2A$ shows that $1A,3B$ and $5B$ are the strongly real $2$-regular classes of $G$. So the real $2$-regular classes $3A$ and $5A$ are weakly real. Now $\varphi_3=\chi_3^*$ and $\chi_3(5B)=1$ is odd. Also $\varphi_{11}$ is the unique real irreducible Brauer character which occurs with odd multiplicity in $\chi_{15}^*$ and $\chi_{15}(3B)=9$ is odd.  So $\Phi_1,\Phi_3$ and $\Phi_{11}$ are the quadratic indecomposable characters of $G$.

Next $\varphi_2$ and $\varphi_{10}$ are the remaining real irreducible Brauer characters of $G$. The symplectic character $\chi_{11}$ occurs with odd multiplicity $1$ in $\Phi_{10}$. This confirms that $\Phi_{10}$ is a non-quadratic principal indecomposable character. There is a unique symplectic irreducible constituent of $\Phi_2$, namely $\chi_{13}$, and this occurs with even multiplicity $2$. So we cannot use its presence to verify that $\Phi_2$ is not quadratic. On the other hand $\phi_2=\chi_2^*$ and the values $\chi_2(1)=22,\chi_2(3B)=4$ and $\chi_2(5B)=2$ on the strongly real $2$-regular classes are all even. This confirms that $\Phi_2$ is not quadratic.
\end{Example}


\begin{Example}
$G=2.Ru$ has $61$ irreducible ${\mathbb C}$-characters, $25$ of which are faithful. Of the faithful characters, $7$ are symplectic and the remaining $18$ are non-real. Of the $36$ characters of $Ru$, two are non-real and the remaining $34$ are orthogonal. Also $G$ has $9$ classes of elements of odd order, all of whom are real. By examining the square of the two classes of non-central involutions, we see that the weakly real $2$-regular classes are $15A,29A$ and $29B$.

Now the symplectic character $\chi_{53}$ has odd multiplicity $1$ in $\Phi_4$. So $\Phi_4$ is not quadratic. Alternatively, $\varphi_4$ is the unique Brauer character which occurs with odd multiplicity in $\chi_{53}^*$ and $\chi_{53}(15A)=1$. Next $\varphi_6=\chi_{49}^*$ and $\chi_{49}(15A)=-1$ is odd. Similarly $\varphi_7=\chi_{50}^*$ and $\chi_{50}(15A)=-1$. This confirms that $\Phi_6$ and $\Phi_7$ are not quadratic.
\end{Example}



We take the opportunity to correct Example 2.12 in \cite{GW93}. This erroneously claims that a certain principal indecomposable module of a group of order $288$ is weakly real. In fact the assertion, on the first line of p268, that a certain form $c$ is $B$-invariant, is false.

\begin{Example}
Let $H$ be the non-abelian group $C_3\rtimes C_4$ of order $12$. Let $\sigma$ be the switching automorphism of $H\times H$: $(x,y)^\sigma=(y,x)$, for all $x,y\in H$. Set $G=(H\times H)\langle\sigma\rangle$. Now $G$ is a $2$-nilpotent group with normal Hall $2'$-subgroup $N\cong C_3\times C_3$. Let $\omega$ be one of the two non-trivial linear characters of $C_3$. Then $G$ has three orbits on $\op{Irr}(N)$, with representatives $1\times1,1\times\omega$ and $\omega\times\omega$, respectively. Set $\Phi_1=(1\times1){\uparrow^G}$, $\Phi_2=(1\times\omega){\uparrow^G}$ and $\Phi_3=(\omega\times\omega){\uparrow^G}$. Let $\Phi_i$ belong to the $2$-block $B_i$ of $G$. Then $B_1,B_2$ and $B_3$ are distinct, real and nilpotent.

The principal $2$-block $B_1$ consists of the $8$ linear characters and $6$ irreducible characters of degree $2$ in $\op{Irr}(G/N)$. Then $\Phi_1$ is strongly real, as it is the trivial principal indecomposable character of $G$.

The block $B_2$ consists of the $8$ irreducible characters in $\op{Irr}(G\mid1\times\omega)$, each of which has degree $4$. Four of these are orthogonal and four are symplectic. Now $1\times\omega$ has stabilizer $H\times C_6$ and extended stablizer $H\times H$ in $G$. As $H\times H$ does not split over $H\times C_6$, it follows that $B_2$ is a weakly real $2$-block of $G$. So $\Phi_2$ is weakly real.

Finally $B_3$ consists of the $5$ irreducible characters in $\op{Irr}(G\mid\omega\times\omega)$. Four of these are orthogonal and of degree $4$. The remaining character is symplectic and of degree $8$. Now the stabilizer of $\omega\times\omega$ in $G$ is $C_6\wr C_2$. As $(b^{-1},b)\sigma$ is an involution in $G$ which inverts $\omega\times\omega$, the extended stablizer of $\omega\times\omega$ splits over its stabilizer. As a consequence $\op{Irr}(G\mid\omega\times\omega)$ is a strongly real $2$-block of $G$. So $\Phi_3$ is strongly real. This means that the corresponding principal indecomposable $kG$ module has a quadratic geometry, contrary to the conclusion of \cite[2.12]{GW93}.

We note that the involution module of $B_3$ can be constructed as follows. We have $\op{Irr}(H\mid\omega)=\{\chi_1,\chi_{-1}\}$, where $\chi_\epsilon$ has degree $2$ and F-S indicator $\epsilon$. Clearly $\op{C}_G(\sigma)=\Delta H\times\langle\sigma\rangle$ in $G$. So ${\mathbb C}_{\op{C}_G(\sigma)}\uparrow^G$ is isomorphic to ${\mathbb C}H$, as a module for ${\mathbb C}H\wr C_2$. In particular the involution module of $B_3$ has ordinary character $\hat \chi_1+\hat\chi_{-1}$, where $\hat\chi_\epsilon$ is an extension of $\chi_\epsilon\times\chi_\epsilon$ from $H\times H$ to $G$. It is easy to check that both of these characters has F-S indicator $+1$.
\end{Example}

In view of the above examples, we venture the following:

\begin{Conjecture}
$\epsilon(\Phi)$ is even if\/ $\Phi$ is a weakly real principal indecomposable character. 
\end{Conjecture}

\section{Odd Cartan Invariants}

We keep our notation $\varphi_i,\Phi_i$ and $g_i$ for the irreducible $2$-Brauer characters, the principal indecomposable characters and the elements of the $2$-regular conjugacy classes of $G$, respectively. List the irreducible characters of $G$ as $\chi_1,\dots,\chi_k$. Then
$$
\chi_i^*=\sum_{j=1}^\ell d_{ij}\varphi_j,\quad
\Phi_i=\sum_{j=1}^k d_{ij}\chi_j,\quad
\Phi_i=\sum_{j=1}^\ell c_{ij}\varphi_j,
$$
where the Cartan invariants  are given by $c_{ij}=\sum_{u=1}^kd_{ui}d_{uj}$.

Recall that $g\in G$ is said to have $2$-defect zero if $\op{C}_G(g)$ has odd order. In particular $g$ has odd order. Now suppose that $g$ is real and of $2$-defect zero. Then a Sylow $2$-subgroup of $\op{C}_G^*(g)$ has order $2$. So $g$ is inverted by an involution, whence $g$ is strongly real in $G$.

Now the $\ell\times\ell$ Cartan matrix $[c_{ij}]$ of $G$ is a symmetric integer matrix whose invariant factors are $|\op{C}_G(g_1)|_2,\dots,|\op{C}_G(g_\ell)|_2$. In particular its rank modulo $2$ coincides with the number of $2$-regular conjugacy classes of $G$ which have $2$-defect zero. So we get one strongly real $2$-Brauer character for each real class of $2$-defect zero. Our final Theorem refines this observation:

\begin{Theorem}
Suppose that $c_{ij}$ is odd, where $\varphi_i$ and $\varphi_j$ are real-valued. Then there exists $w$ such that $c_{iw}$ is odd and $\varphi_w$ is strongly real. 
\begin{proof}
We compute
$$
\begin{aligned}
\sum_{u=1}^r\Phi_i(g_u)\frac{\Phi_j(g_u^{-1})}{|\op{C}_G(g_u)|}
&\equiv\sum_{u=1}^r\sum_{v=1}^\ell c_{iv}\varphi_v(g_u)\frac{\Phi_j(g_u^{-1})}{|\op{C}_G(g_u)|},\quad\mbox{definition of $c_{iv}$}\\
&\equiv\sum_{u=1}^r\sum_{v=1}^rc_{iv}\varphi_v(g_u)\frac{\Phi_j(g_u^{-1})}{|\op{C}_G(g_u)|},\quad\mbox{as $c_{iv}\varphi_v(g_u)=c_{i\ov v}\ov\varphi_v(g_u)$}\\
&\equiv\sum_{v=1}^rc_{iv}\sum_{u=1}^\ell\varphi_v(g_u)\frac{\Phi_j(g_u^{-1})}{|\op{C}_G(g_u)|},\quad\mbox{as $\varphi_v(g_u)\Phi_j(g_u^{-1})=\varphi_v(g_u^{-1})\Phi_j(g_u)$}\\
&\equiv c_{ij}\mod J,\quad\mbox{by the second orthogonality relation}.
\end{aligned}
$$
It follows that $\Phi_i(g_u)\frac{\Phi_j(g_u^{-1})}{|\op{C}_G(g_u)|}\not\in J$, for some real $g_u$. But $|\op{C}_G(g_u)|_2\mid\Phi_i(g_u)$. So $g_u$ must have $2$-defect zero, and hence $g_u$ is strongly real. Now $\Phi_i(g_u)=\sum_{w=1}^\ell c_{iw}\varphi_w(g_u)\equiv\sum_{w=1}^rc_{iw}\varphi_w(g_u)$. So there exists $w$ such that $\varphi_v=\ov\varphi_v$, $c_{iw}$ is odd and $\varphi_w(g_u)$ is coprime to $2$. Then $\varphi_w$ is strongly real, according to Lemma \ref{L:strongly_real}. This completes the proof.

{\bf Remarks:} Let $t$ be an involution inverting a primitive idempotent corresponding to $\Phi_w$. Then $t$ inverts $g_u$, by Lemma \ref{L:inverts}. So $t$ is uniquely determined up to conjugacy, as $g_u$ has defect $0$. Now $c_{iw}=\sum_{x=1}^kd_{ix}d_{wx}$. So there exists $x$ such that $\chi_x=\ov\chi_x$ and $d_{ix}$ and $d_{wx}$ are odd. As $\varphi_w$ is strongly real, $\chi_x$ is an orthogonal irreducible character of $G$.
\end{proof}
\end{Theorem}

\section{bibliography}

\end{document}